\documentclass[12pt]{article}
\usepackage{latexsym, amsmath, amsfonts, amsthm}

\newcommand{\R}{{\mathbb R}}                  
\newcommand{\N}{{\mathbb N}}                      
\newcommand{\Z}{{\mathbb Z}}                     
\newcommand{\C}{{\mathbb C}}

\newcommand{\dis}{\displaystyle}
\newcommand{\E}{{\mathbb E}}

\newcommand{\B}{{\cal B}}
\newcommand{\lan}{\langle}
\newcommand{\ran}{\rangle}

\newcommand{\diam}{\mathop{\rm diam}}
\newcommand{\spa}{\mathop{\rm span}}

\newcommand{\al}{\alpha}
\newcommand{\sig}{\sigma}
\newcommand{\eps}{\varepsilon}
\newcommand{\vphi}{\varphi}
\newcommand{\la}{\lambda}
\newcommand{\de}{\delta}
\newcommand{\Ee}{{\cal E}}

\newcommand{\pr}[2]{\langle {#1} , {#2} \rangle}
\newcommand{\norm}[1]{\left \| #1 \right \|}

\newtheorem{theorem}{Theorem}
\newtheorem*{theorem*}{Theorem}

\newtheorem{lemma}{Lemma}
\newtheorem{proposition}{Proposition}
\newtheorem{definition}{Definition}
\newenvironment{thm}[2][]
               {\pagebreak[3] \medskip {\bf \noindent Theorem}
                               {#1} {\it #2}
                             }{}
\title{Majorizing measures and proportional subsets of bounded
orthonormal systems}
\author{Olivier GU\'{E}DON \and Shahar MENDELSON${}^{1}$ \and Alain PAJOR
       \and
Nicole TOMCZAK-JAEGERMANN${}^{2}$ }

\newcommand\address{\noindent\leavevmode%
\noindent
 O. Gu\'{e}don,\\
 Universit\'e Pierre et Marie Curie,
  Paris 6, \\
 Institut de   Math\'ematiques de Jussieu,\\
4 place Jussieu, 75005 Paris, France\\
\texttt{\small%
e-mail:  guedon@math.jussieu.fr}

\medskip
\noindent
S. Mendelson,\\
CMA, The Australian National University, \\
Canberra, ACT 0200,  Australia, and \\
Department of Mathematics, Technion, \\
Haifa, 32000, Israel \\
\texttt{\small%
e-mail: shahar.mendelson@anu.edu.au}

\medskip
\noindent
A. Pajor, \\
Universit\'{e} Paris-Est\\
\'{E}quipe d'Analyse et Math\'{e}matiques Appliqu\'ees, \\
5, boulevard Descartes,
Champs sur Marne,\\
77454 Marne-la-Vall\'{e}e,  Cedex 2, France\\
\texttt{\small%
e-mail: Alain.Pajor@univ-mlv.fr }

\medskip
\noindent
N. Tomczak-Jaegermann, \\
Dept.~of Math.~and Stat.~Sciences,\\
University of Alberta, \\
Edmonton, Alberta, Canada, T6G 2G1.\\
\texttt{\small%
e-mail:    nicole.tomczak@ualberta.ca}
}

\date{}

\begin{document}
\maketitle

\maketitle \footnotetext{{\it 2000 MSC-classification:}
    46B07,  46B09, 42A05, 42A61}
\footnotetext[1]{Partially supported by an Australian Research
Council Discovery grant and an Israel Sciences foundation grant 666/06.}
\footnotetext[2]{This author holds the
Canada Research Chair in Geometric Analysis.}

\begin{abstract}
In this article we prove that for any orthonormal system
$(\vphi_j)_{j=1}^n \subset L_2$ that is bounded in $L_{\infty}$,
and any $1 < k <n$, there exists a subset $I$ of cardinality greater than
$n-k$ such that on
$\spa\{\vphi_i\}_{i \in I}$, the $L_1$ norm and the $L_2$ norm are
equivalent up to a factor $\mu (\log \mu)^{5/2}$, where $\mu =
\sqrt{n/k} \sqrt{\log k}$. The proof is based on a new estimate of
the supremum of an empirical process on the unit ball
of a Banach space with a good modulus of convexity, via the use of
majorizing measures.
\end{abstract}

\section{Introduction.}

We study some natural empirical processes determined by uniformly
convex Banach spaces with modulus of convexity of power type 2.
Results of this kind were extensively studied in a Hilbertian
setting, and became an important tool for investigations, for
example, of the behaviour of various random sets of vectors (as in
\cite{Rudelson, Rudelson1, GueRud, RudelsonVershynin, CT1, GMPT}).
We then apply these results to address a problem of selecting a
subset of a bounded orthonormal system (for example, a set of
characters), in the spirit of a result of Bourgain (see
\cite{Tal_proport}) and of Talagrand \cite{Tal_proport}, that also
has applications to sparse reconstruction (\cite{CT1,
RudelsonVershynin, GMPT}). A particular case of this result,
formulated in Theorem \ref{Lambda_2}, is a so called Kashin's
splitting of a set of $2k$ orthonormal vectors in $L_2$ that is
bounded in $L_{\infty}$.
\begin{theorem*}
There exist two positive constants $c$ and $C$ such that for
any even integer $n$ and
any orthonormal system $(\vphi_j)_{j=1}^n$  in $ L_2$ with
$\|\vphi_j\|_{L_\infty} \le L$ for $1 \le j \le n$,
we can find  a subset $I \subset \{1, \dots, n\}$
with $n/2 - c \sqrt n \le |I| \le  n/2 + c \sqrt n$
such that for every $a = (a_i) \in \C^n$,
$$
\left\|\, \sum_{i\in I} a_i \vphi_i\, \right\|_{L_2}
\le C \ L\,  \sqrt{\log n} (\log \log n)^{5/2}
\left\|\, \sum_{i\in I} a_i \vphi_i\, \right\|_{L_1}
$$
and
$$
\left\|\, \sum_{i\notin I} a_i \vphi_i\, \right\|_{L_2}
\le C \ L\,  \sqrt{\log n} (\log \log n)^{5/2}
\left\|\, \sum_{i\notin I} a_i \vphi_i\, \right\|_{L_1}.
$$
\end{theorem*}
This result strengthens  Theorem 2.4 of \cite{GMPT} and is almost
optimal.  We have been
told by J. Bourgain that the term $\sqrt{\log n}$ is necessary.
We would like to thank
L. Rodriguez-Piazza for showing us the details of the proof of this optimality in the case of the Walsh system and
for allowing us to present the argument at the end of this paper.
The  technical proof of  Theorem \ref{majorizing} about empirical processes
will be presented in the first part. It is based
on a construction of majorizing measures developped in \cite{Rudelson} and \cite{GueRud} and the main
new ideas we use are some
observations on packing and covering numbers.

\section{Maximal deviation of the empirical moment.}
We begin this section with some definitions and notation.
If $E$ is a normed space we
denote by $E^*$ the dual space to $E$ and the dual norm is denoted
by $\| \ \|_*$.
The modulus of convexity of $E$ is defined for any $\eps \in (0,2)$
by
$$
  \delta_E (\eps) = \inf \left \{1 - \norm{\frac{x+y}{2}}, \|x\| = 1,
  \|y\| = 1, \|x-y\| > \eps \right \}.
$$
We say that $E$ has modulus of convexity of power type $2$ if there
is a constant $c$ such that $\delta_E(\eps) \ge c \eps^2$ for every
$\eps \in (0,2)$. It is well-known (see e.g., \cite{Pisier},
Proposition 2.4)
that this property is equivalent to the fact that the inequality
\begin{equation}
    \label{eq:powertype2}
\norm{ \frac{x+y}{2}}^2 + \la^{-2} \norm{ \frac{x-y}{2}}^2 \le
\frac{1}{2} (\|x\|^2 + \|y\|^2)
\end{equation}
holds for all $x,y \in E$ (where $\la>0$ is a constant depending
only on $c$). If $(\ref{eq:powertype2})$ is satisfied then we
say that $E$ has modulus of convexity of power type $2$ with
constant $\la$ (in such a case it  is clear that $\delta_E(\eps)
\ge \eps^2/2\la^2$).

The notions of type (and cotype) of a Banach space were studied
extensively during the 70's (see, for example, \cite{MaureyPisier}
and the survey \cite{Maurey}).
A Banach space $E$ has type $p$ if there is a constant $C$ such that
for every $N \in \N$ and every $x_1, \ldots, x_N$
\begin{equation} \label{eq:type}
\E\| \sum_{i=1}^N g_i x_i \|
\le
C \left( \sum_{i=1}^N \|x_i\|^p \right)^{1/p}
\end{equation}
where $g_1, \ldots, g_N$ are standard independent gaussian variables
(that is $g_i \sim {\cal N} (0,1)$). The smallest constant $C$ for
which \eqref{eq:type} holds is called the type $p$ constant of $E$
and is denoted by $T_p(E)$.

Moreover, it is well known
that if $E$ has modulus of convexity of power type 2
then the dual space $E^*$ has also modulus of smoothness of power
type 2, and therefore, $E^*$ has type 2 (see, for example, Theorem 1.e.16
of \cite{LinTza}).
%
%
\subsection{Results on empirical processes.}
Our first theorem generalizes
a result of Rudelson \cite{Rudelson1}, proved in a Hilbertian
setting, to the case of a Banach space with modulus of convexity of
power type 2. In fact, in Theorem \ref{majorizing} we solve a
question left open in \cite{GueRud} by removing the condition on the
distance of $E$ to an Euclidean space of the same dimension.
\begin{theorem}\label{majorizing}
There exists an absolute constant $C$ for which the following holds.
Let $E$ be a Banach space with modulus of convexity of power type
$2$ with constant $\la$. Then, for every  vectors
$X_{1}, \ldots, X_{m}$ in $E^*$,
\begin{align} \label{eq:power2}
\E & \sup_{y \in B_E} \left|
 \sum_{j=1}^m \eps_{j} |\lan X_j, y \ran|^2
 \right| \le
\\
\nonumber & C \, \dis \la^4 \, T_2(E^*) \sqrt{\log m}
 \max_{1 \le j \le m} \|X_{j}\|_{*}
 \sup_{y \in B_E} \left(\sum_{j=1}^m |\lan X_j, y
 \ran|^2 \right)^{1/2},
\end{align}
 where the expectation is taken over the i.i.d. Bernoulli random
variables $(\eps_{j})_{1 \le j \le m}$ and $T_2(E^*)$ is the type
$2$ constant of $E^*$.
\end{theorem}
 %
%
%
%
The proof of Theorem \ref{majorizing} uses the same construction of
a majorizing measure used in \cite{Rudelson} and in \cite{GueRud},
and this construction was inspired by the work of Talagrand in
\cite{Talagrand}. The improvement in Theorem \ref{majorizing}
compared with \cite{GueRud} comes from entropy estimates that will
be presented in Section \ref{entropy}. We will only sketch the
construction of the majorizing measure in Section \ref{construction}
and explain how the argument from \cite{GueRud} may be adapted using
the new estimates.

As it has been proved in \cite{Rudelson1}, one can apply a
symmetrization argument due to Gin\'e and Zinn \cite{GineZinn}
combined with Theorem \ref{majorizing} and deduce the following
result.
\begin{theorem}
    \label{thrandom}
There exists an absolute constant $C$ for which the following holds.
Let $E$ be a Banach space with modulus of convexity of power type
$2$ with constant $\la$. Let $X$ be a random vector in $E$ and set
$X_1, \ldots, X_m$ to be independent copies of $X$. If
$$
A = C\la^4 \, T_2(E^*) \sqrt{\frac{\log m}{m}} \left(\E \max_{1 \le
j \le m} \|X_j\|_{*}^2\right)^{1/2} \hbox{and } \sig^2= \sup_{y \in
B_E} \E|\lan X, y \ran|^2
$$
then
$$
 \E
 \sup_{y \in B_E}
 \left| \frac{1}{m} \sum_{j=1}^m |\lan X_j, y \ran|^2 - \E |\lan X, y \ran|^2
 \right| \le A^2 +  \sig \, A.
$$
\end{theorem}
We omit further details and refer
the reader to one of the articles
\cite{Rudelson1,GueRud,RudelsonVershynin,GMPT}.


\subsection{Covering and packing numbers.}
\label{entropy}
The new ingredient of our proof, in comparison with \cite{GueRud},
are estimates on packing and covering numbers which we shall now
discuss.

\begin{definition} Let $T$ and $B$ be
symmetric convex bodies in a Banach space $E$. Define $N(T,B)$ to be
the minimal number of translates of $B$ needed to cover $T$ i.e.
$$
N(T,B) = \inf \{N, \exists \{x_{1}, \ldots, x_N\} \subset T \hbox{
such that } T \subset \cup_{j=1}^N x_j + B \}.
$$
We denote by $M(T,B)$ the maximal number of disjoint translates of
$B$ by elements of $T$ i.e.
$$
M(T,B) = \sup \{N, \exists \{x_{1}, \ldots, x_N\} \subset T \hbox{
such that } x_j - x_k \notin B \hbox{ for } j \ne k\}.
$$
\end{definition}
It is well known that these two quantities are related by
\begin{equation}
    N(T,B) \le M(T,B) \le N(T, B/2).
    \label{eq:entropypacking}
\end{equation}
Note that if $E$ is a normed space, $B_E$ is its unit ball and
$B=\varepsilon B_E$, then $M(T,\varepsilon B_E)$ is the cardinality
of a maximal $\varepsilon$-separated subset of $T$ with respect to
the norm in $E$. Also, if $S:X \to Y$ is an operator, we sometimes
denote $M(S B_X,\varepsilon B_Y)$ by $M(S:X \to Y, \varepsilon)$,
and use a similar notation for the covering numbers $N$.

The following is a combination of results that appeared in \cite{BPST}
and in \cite{C} (see Lemma 3.3 in \cite{Talagrand}).
We repeat its proof since
we need precise estimates on the dependence on the modulus of
convexity of the space $E$.
\begin{lemma}\label{cov-ellinfini}
Let $E$ be a  Banach space of modulus of convexity of power type 2
with constant $\la$
and denote by $T_2(E^*)$ the type 2 constant of $E^*$.
\\
Let $X_1, \ldots, X_m$ be vectors in $E^*$ such that for every $1
\le j \le m$, $\|X_j\|_{*} \le L$, and define for every $y \in E$,
$$
\|y\|_{\infty, m} = \max_{1 \le \ell \le m} | \lan X_{\ell}, y \ran |.
$$
If $\eps > 0$ and $x_1, \ldots, x_N \in B_E$ are $\eps$-separated
with respect to $\| \cdot \|_{\infty,m}$ then
$$
\eps \sqrt{ \log N} \le C \la^{2} \ T_2(E^*) \ L \, \sqrt{\log m}
$$
where $C$ is an absolute constant.
\end{lemma}
{\bf Proof.} Let $(e_i)_{i=1}^m$ be the standard basis in $\ell_1^m$
and define $S: \ell_1^m \to E^*$ by $S e_{\ell} = X_{\ell}$ for
every $\ell=1, \ldots, m$. Hence, for every $\ell \ne \ell'$,
$|S^*(x_{\ell} - x_{\ell'})|_{\infty} \ge \eps$ where
$|\cdot|_{\infty}$ denotes the norm in $\ell_{\infty}^m$. Thus,
$$
N \le  M(S^*:E \to \ell_\infty^m, \eps).
$$
Since
$\|S\| = \max_{1 \le \ell \le N} \| X_{\ell}\|_* \le L$ and $E^*$ has
type 2, Proposition  1 in \cite{C} yields that
$$
\eps \sqrt{\log N(S: \ell_1^m \to E^*, \eps)} \le C \,  L \,
T_2(E^*) \sqrt{\log m}
$$
where $C$ is an absolute constant. By $(\ref{eq:entropypacking})$, a
similar estimate holds for $M(S: \ell_1^m \to E^*, \eps)$, where $C$
is replaced by a new constant, also denoted by $C$. Define the
function
$$
f(\eps) = \eps \sqrt{\log M(S: \ell_1^m \to E^*, \eps)}
$$
and observe that $f$ is bounded by $C \,  L \, T_2(E^*) \sqrt{\log m}$.

Since $E$ has modulus of convexity of power type 2 with constant
$\la$, Proposition 2 in \cite{BPST} states that for every $\theta
\ge \eps$,
$$
M(S^* : E \to
\ell_{\infty}^m, \eps)
\le
M(S^* : E \to
\ell_{\infty}^m, \theta)
M(S: \ell_1^m \to E^*, c \theta \delta_E(\eps/ \theta))
$$
and recall that $\de_E(\eps/ \theta) \ge \eps^2 / 8 \la^2 \theta^2$.
Therefore, if
$$
h(\eps) = \eps \sqrt{ \log M(S^* : E \to
\ell_{\infty}^m, \eps)} ,
$$
then the previous inequality implies that
$$
h(\eps) \le \frac{\eps}{\theta} \,  h(\theta) + c^\prime \, \left(
\frac{\la^2 \theta}{\eps} \right)  f(\frac{c\eps^2}{\theta \la^2})
$$
where $c,c^\prime$ are absolute constants.

Choosing $\theta = 2 \eps$, it follows that $h$ is bounded by
$$
C^\prime \, \la^2 \,  L \, T_2(E^*) \sqrt{\log m}
$$
proving the announced result. $\hfill \Box$

The next lemma is a simple application of Sudakov's  inequality
\cite{Sudakov}.
\begin{lemma}\label{cov-newnorm}
There exists an absolute constant $C$ for which the following holds.
Let $E$ be a Banach space such that $E^*$ has type 2. Assume that
$X_1, \ldots, X_m$ are vectors in $E^*$ with $\sup_{1 \leq j \leq m}
\|X_j\|_{*} \le L$. Let $\al_1, \ldots, \al_m$ such that
$\sum_{\ell = 1}^m \al_{\ell}^2 \le 1$ and
for every $y \in E$ set
$$
|y|_{{\cal E}}^2 = \sum_{\ell = 1}^m
| \langle X_{\ell}, y  \rangle |^2
\al_{\ell}^2.
$$
If $\eps >0$ and $x_1, \ldots, x_N \in B_E$ are $\eps$-separated
with respect to $|\cdot |_{\cal E}$, then
$$
\eps \sqrt{ \log N} \le C T_2(E^*) L.
$$
\end{lemma}
{\bf Proof.} Let $\Ee_1$ be the ellipsoid in $\R^m$
consisting of all $y \in \R^m$ such that
$|y|_{\Ee_1} := \left(\sum_{\ell = 1}^m \al_{\ell}^2 \lan y,
e_{\ell}\ran^2 \right)^{1/2} \le 1$. Set $H = (\R^m, |\cdot|_{\Ee_1})$
and define the operator $S: H^* \to E^*$
by $Se_i = X_i$, where $(e_i)_{i=1}^m$ is the standard basis in
$\R^m$. Note that for every $x_i,x_j$,
$$
|S^*(x_i - x_{j}) |_{\Ee_1}^2 = \sum_{\ell = 1}^m \al_{\ell}^2 \lan
S^*(x_i - x_{j}), e_{\ell} \ran^2 = |x_{i} - x_{j}|_{\Ee}^2,
$$
and thus the points $\{S^*(x_{1}), \ldots, S^*(x_N)\}$ are
$\eps$-separated in $| \cdot|_{\Ee_1}$ and belong to $S^*(B_E)$. By
$(\ref{eq:entropypacking})$, $N \le N(S^*(B_E), (\eps / 2) \Ee_1)$. On the
other hand, $\Ee_1 = T^{-1}(B_2^m)$ where $T: \R^m \to \R^m$ is the
diagonal operator $T e_{\ell} = \al_{\ell} \, e_{\ell}$. Applying
Sudakov's inequality \cite{Sudakov},
\begin{equation}
\label{N2}
\eps \sqrt{\log N} \le
\eps \sqrt{\log N(S^*(B_E), (\eps / 2)  T^{-1}(B_2^m))} \le C
\E \sup_{z \in T(S^*(B_E))} \lan G, z \ran
\end{equation}
where $G$ is a canonical Gaussian vector in $\ell_2^m$ and $C$ is an
absolute constant. Moreover
\begin{eqnarray*}
\E \sup_{z \in T(S^*(B_E))} \lan G, z \ran & = & \E
\sup_{y \in B_E} \lan T G, S^* y \ran
\\
 =
\E \| \sum_{\ell = 1}^m \al_{\ell} g_{\ell} X_{\ell} \|_{*}
& \le &
T_2(E^*) \left( \sum_{\ell = 1}^m \al_{\ell}^2 \|X_{\ell} \|_{*}^2 \right)^{1/2}
\le
T_2(E^*) L
\end{eqnarray*}
where we have used the type 2 inequality for $E^*$ and the fact that
for every $1 \leq \ell \leq m$, $\|X_{\ell} \|_{*} \le L$. The
result now follows from  $(\ref{N2})$. $\hfill \Box$

\subsection{Construction of a majorizing measure.}
\label{construction} The construction that we present here is the
same one that was presented in \cite{Rudelson} and in \cite{GueRud}.
Let $X_1, \ldots, X_m$ be $m$ fixed vectors in $E^*$ and define the
random process $\{V_y : y \in B_E\}$ by
$$
V_y = \sum_{j=1}^m \eps_j | \langle X_j,y \rangle |^2,
$$
where $\eps_j$ are independent symmetric Bernoulli random variables.

Our aim is to show that when $E$ has modulus of convexity of power
type~2,
\begin{eqnarray}\label{sup}
\E \sup_{y \in B_E} |V_y| \le C\,  \la^4 \, \, T_2(E^*) \, \max_{1
\leq j \leq m} \|X_j\|_* \, \sqrt{\log m} \left(\sup_{y \in B_E}
\sum_{i=1}^m |\pr{X_j}{y}|^2\right)^{1/2}
\end{eqnarray}
for a suitable absolute constant $C$.

It is known that the process $\{V_y : y \in B_E\}$ satisfies a
sub-Gaussian tail estimate, namely, that for every $y, \overline y
\in E$ and any $t > 0$,
$$
P(|V_y - V_{\overline y}| \ge t)
\le
2
\exp\left(- \frac{c t^2}{\tilde{d}^2(y,\overline y)} \right)
$$
where
$$
\tilde{d}^2(y,\overline y) = \sum_{j=1}^m  \left( | \langle X_j,y
\rangle |^2 - | \langle X_j, \overline y \rangle |^2 \right)^2
$$
and $c$ is an absolute constant.

It will be preferable to consider the following quasi-metric
$$
d^2(y,\overline y) = \sum_{j=1}^m  | \langle X_j,y-\overline y
\rangle |^{2} \left( | \langle X_j,y \rangle |^2 + | \langle X_j,
\overline y \rangle |^2 \right)
$$
and the quasi-norm $\| \cdot \|_{\infty, m}$ endowed on $E$ by
$$
\|x\|_{\infty, m} = \max_{1 \le j \le m} |\lan X_j,x\ran|.
$$
The proof of inequality $(\ref{sup})$ is based on the majorizing
measure theory of Talagrand \cite{Tal:genericchaining}. The
following theorem is a combination of Proposition 2.3, Theorem 4.1
and Proposition 4.4 of \cite{Tal:genericchaining}.

\begin{thm}[\cite{Tal:genericchaining}]
{\it Assume that the process $\{V_y : y \in B_E\}$ is subgaussian
with respect to a metric $d$. Let $r \ge 2$ and $k_0 \in \Z$ be the
largest integer such that
$r^{-k_0}$ is greater than the radius of $B_E$ with respect to the
metric $d$. For
every $k \geq k_0$ let $\phi_k : B_E \to \R^+$ be a family of maps
satisfying the following assumption: there exists $A>0$ such that
for any point $x \in B_E$, any $k \ge k_0$ and any $N \in \N$
\\
$ (H) \left\{
\begin{array}{l}
for \ any \ points \ x_1, \ldots, x_N \in \B_{r^{-k}}(x) \ with \
d(x_i, x_j) \ge r^{-k-1}, i \ne j
\\
we \ have \ \dis \max_{i=1, \ldots, N} \phi_{k+2}(x_i) \ge \phi_k(x)
+ \, \frac{1}{A} \ r^{-k} \sqrt{\log N}.
\end{array}
\right. $
\\
Then
$$
  \E \sup_{y \in B_E} |V_y - V_0|
  \le c \ A \cdot \sup_{k\ge k_0, x\in K} \phi_k(x).
$$
}
\end{thm}

The construction requires certain properties of the quasi metric $d$
and the quasi norm $\| \ \|_{\infty,m}$. We refer to Propositions 1
and 2 in \cite{GueRud} for precise properties of these metrics and
list the ones we require in the following lemma.
\begin{lemma} \label{lemma:prop}
For every $y, \overline y \in \R^n$ and every $u \in B_E$,
\begin{eqnarray}
\label{ineg1} & \tilde{d}(y, \overline y) \le 2 \, d(y,\overline y),
\\
\label{ineg2} & d(y,\overline y) \le  \sqrt 2 \|y-\overline
y\|_{\infty, m} \ \sqrt M,
\\
\label{ineg3} & \|y - \overline y\|_{\infty, m}  \le \max\limits_{1 \le j
\le m}\|X_j\|_* \ \|y-\overline y\|,
\\
\label{ineg4} & d^2(z, \overline z) \le 8 \left( |z-\overline
z|_{{\cal E}_u}^2 +
M \,  \|z-\overline z\|_{\infty, m}^2
(\|z-u\|^{2} +\|\overline z -u \|^{2} )\right),
\end{eqnarray}
where
$$
|x|_{{\cal E}_u}^2=\sum_{i=1}^m
{\pr{X_i}{x}}^2{\pr{X_i}{u}}^2 \hbox{ and }
M = \sup_{y \in B_E} \sum_{j=1}^m | \langle X_j,y \rangle |^2.
$$
Moreover, for every $x \in E$ and $\rho>0$, the ball (with respect
to the quasi-metric $d$) centered in $x$ and with radius $\rho$,
denoted by $\B_x(\rho)$, is convex.
\end{lemma}
Note that by combining \eqref{ineg2} and \eqref{ineg3} it follows
that for every $\rho>0$ and every $x \in B_E$, $\inf_{y \in
\B_\rho(x)} \|y\|$ is attained.

%
%
\bigskip
{\bf Proof of Theorem \ref{majorizing}.} Since there is only a
finite number of points $X_1, \ldots, X_m$ then by passing to a
quotient of $E$ we can assume that $E$ is a finite dimensional
space. We will denote its dimension by $n$ and obviously one can
assume that $m \ge n$. Also, recall that if $E$ is a Banach space
with modulus of convexity of power type $2$ with constant $\la$ then
every quotient of $E$ satisfies that property with a constant
smaller than $\la$ (see, e.g., \cite{LinTza}).

By the homogeneity of the statement we can assume that
\begin{equation}
    \label{homo}
\sup_{y \in B_E} \sum_{j=1}^m | \langle X_j,y \rangle |^2 = 1
\end{equation}
and by inequality (\ref{ineg1}), $V_y$ is a sub-Gaussian process
with the quasi-metric $2  d$.

Therefore, if we denote
$$
L = \max_{1\le j \le m} \|X_j\|_*
$$
our aim is to show that
\begin{eqnarray*}
\E \sup_{y \in B_E} |V_y| \le C\,  \la^4 \, \, T_2(E^*) \, L \,
\sqrt{\log m}
\end{eqnarray*}
for an absolute constant $C$.

By inequality $(\ref{ineg2})$, the diameter of $B_E$ with respect to
the metric $d$ is bounded by $2\sqrt 2 L$. Let $r$ be a fixed number
chosen large enough, set $k_0$ to be the largest integer such that
$r^{-k_0} \ge 2\sqrt 2 L$ and put $k_1$ to be the smallest integer
such that $r^{-k_1} \le L / \sqrt n$, where $n$ is the dimension of
$E$. We shall use the same definition of the functionals $\phi_k :
B_E \to \R^+$ as in \cite{Rudelson} and \cite{GueRud}, namely:
$$
\left\{
\begin{array}{l}
\dis \forall k \ge k_1 +1, \phi_k(x) = 1 + \frac{1}{2 \log r} +
\frac{\sqrt n}{L \, \sqrt{\log m}} \sum_{\ell=k_1}^k r^{-\ell}
\sqrt{\log(1 + 4 L  r^{\ell})}
\\
\dis \forall k_0 \le k \le k_1,
\phi_k(x) =
\min\{ \|y\|^2, y \in \B_{8r^{-k}}(x) \} +
\frac{k-k_0}{\log m}.
\end{array}
\right.
$$
It is easy to verify using definitions of $k_0$ and $k_1$ that
$$
\sup_{x \in B_E, k\ge k_0} \phi_k(x) \le c,
$$
where $c$ is an absolute constant.

It remains to prove that our functionals satisfy condition $(H)$ for
$$
A = C \, \la^4 \, L \, \sqrt{\log m},
$$
and that will conclude the proof of Theorem \ref{majorizing}.
$\hfill \Box$
\\
{\bf Proof of condition $(H)$.} Fix integers $N$ and $k$, let $x \in
B_E$ and $x_1, \ldots, x_N \in \B_{r^{-k}}(x)$ for which $d(x_i,
x_j) \ge r^{-k-1}$.

For $k \ge k_1-1$, we always have
$$
\phi_{k+2}(x_i) - \phi_k(x)
\ge
\frac{\sqrt {n \log(1 + 4 L r^{k+2})}}
{L \sqrt{\log m}} \ r^{-k-2}.
$$
Since the points $x_1, \ldots, x_N$ are well separated with respect
to the metric $d$, then by $(\ref{ineg2})$ and $(\ref{ineg3})$,
\[
 \|x_i - x_j\| \ge r^{-k-1} / L \sqrt 2.
\]
By a classical volumetric estimate (see, for example
\cite{PisierBook}), $ N(B_E, t B_E) \le \left ( 1+ 2/t \right )^n $
where $n$ is the dimension of $E$. Therefore,
$$
\sqrt{\log N}  \le \sqrt{n \log(1 + 2 \sqrt{2} L r^{k+1})},
$$
which proves the desired inequality.

The case $k_0 \le k \le k_1-2$ is more delicate and the main
ingredients in this part are the entropy estimates proved in part
\ref{entropy}.

For $j = 1, \ldots, N$  denote by $z_j \in B_E$ points at which
$\min \{ \|y\|^2, y \in \B_{8 r^{-k-2}}(x_j) \} $ is attained and
set $u \in B_E$ to be a point at which $\min \{ \|y\|^2, y \in
\B_{8r^{-k}}(x) \}$ is attained. Put
$$
\theta = \max_{j} \|z_j\|^2 - \|u\|^2,
$$
and then $ \max_{j} \phi_{k+2} (x_j) - \phi_k(x) = \theta +
\frac{2}{\log m} $. We shall prove that
\begin{eqnarray}\label{wanted}
\theta + \frac{2}{\log m}
\ge
r^{-k} \sqrt{\log N} / A.
\end{eqnarray}
Following \cite{GueRud}, it is evident from the properties of ${\cal
B}_{\rho}$ (see proposition 1 and 2 in  \cite{GueRud} and Lemma
\ref{lemma:prop}) that for any $i \ne j$, $d(z_i, z_j) \ge c
r^{-k-1}$ and that $d(x, z_j) \le 8 r^{-k}$. It implies that $(z_j +
u) /2 \in \B_{8r^{-k}}(x)$ and by the definition of $u$, $\|u\| \le
\|z_j + u\| /2$. Since $B_E$ has modulus of convexity of power type
$2$, then for all $j=1, \ldots, N$ $\|z_j - u \| \le \la \sqrt{2
\theta}. $ Thus, the set
$$
U = u +  \la \sqrt{2 \theta} B_E
$$
contains all the $z_j$'s.

Fix an absolute constant $\tilde{c}$ to be named later, set
$$
\delta = \tilde{c} \la^{-1} r^{-k} \theta^{-1/2}
$$
and let $S$ be the maximal number of points in $U$ that are $\delta$
separated in $\|\cdot\|_{\infty, m}$. By Lemma \ref{cov-ellinfini},
$$
\delta \sqrt{\log S} \le
C \, \la^2 \, L \sqrt{\log m} \ \la \ \sqrt \theta
$$
where $C$ is an absolute constant. We consider now two cases.

First, assume that $S \ge \sqrt N$. Then by the previous estimate
and the definition of $\delta$,
$$
\sqrt{\log N}
\le
C \, L \sqrt{\log m}  \ \la^3 \ \sqrt \theta  / \, \delta
\le C \, \theta \
\ r^k \, \la^4  \, L \sqrt{\log m}
$$
which proves  $(\ref{wanted})$.

The second case is when $S\le \sqrt N$. Since $S$ is the maximal
number of points in $U$ that are $\delta$ separated with respect to
$\|\cdot\|_{\infty, m}$, $U$ is covered by $S$ balls of diameter
smaller than $\delta$ in $\|\cdot\|_{\infty, m}$. Thus, there exists
a subset $J$ of $\{1, \ldots, N\}$ with cardinality $|J| \ge \sqrt
N$ such that
$$
\forall i,j \in J, \ \ \|z_i - z_j\|_{\infty, m} \le \delta.
$$
Recall that for any $y \in E$,
$$
|y|_{{\cal E}_u}^2 = \sum_{\ell = 1}^m
| \langle X_{\ell}, y  \rangle |^2
| \langle X_{\ell}, u \rangle |^{2}
=
\sum_{\ell = 1}^m
| \langle X_{\ell}, y  \rangle |^2
\al_{\ell}^{2}.
$$
It is evident from $(\ref{homo})$ that $\sum \al_{\ell}^{2} \le 1$
and from $(\ref{ineg4})$ that for every $z, \overline z , u \in E$,
$$
d^2(z, \overline z) \le
8 \big( |z-\overline z|_{{\cal E}_u}^2 +
 \|z-\overline z\|_{\infty, m}^2
(\|z-u\|^{2}
+\|\overline z -u \|^{2} )\big).
$$
Since $d(z_i, z_j) \ge  c r^{-k-1}$,
$\| z_i - u \| \le 2 \la \sqrt \theta$ and
$\|z_i - z_j\|_{\infty, m} \le \delta$ for any $i, j \in J$,
we can define $\tilde{c}$ small enough
such that for all  $i \ne j \in
J$,
$$
|z_i - z_j|_{{\cal E}_u} \ge c r^{-k-1}.
$$
Thus, there are $|J|$ points in $u +\la \sqrt{2 \theta} B_E$ that
are $c r^{k-1}$-separated for $|\cdot|_{{\cal E}_u}$. Using Lemma
\ref{cov-newnorm},
$$
\sqrt{\log |J|} \le C T_2(E^*) L \ r^k  \ \la \sqrt \theta,
$$
and a simple application of the arithmetic-geometric means
inequality proves that
$$
\sqrt{\log N} \le C T_2(E^*) L \ r^k \ \la \sqrt{\log m}
\left(\theta + \frac{2}{\log m}\right),
$$
completing the proof of $(\ref{wanted})$ (because $\la \ge 1$).
$\hfill \Box$

%
\section{Selecting an arbitrary proportion of a bounded orthonormal system.}
In this section, we will prove the
\begin{theorem}
    \label{Lambda_2}
There exist two positive constants $c, C$ such that for
any orthonormal system $(\vphi_j)_{j=1}^n$ in $ L_2$ with
$\|\vphi_j\|_{L_\infty} \le L$ for $1 \le j \le n$, the following
holds.

1) For any
   $1 < k < n$ there exists a subset $I \subset \{1, \dots, n\}$
   with $|I| \ge n - k $ such that for every $a = (a_i) \in \C^n$,
\begin{equation}
  \label{lambda_2_ineq}
\left\|\, \sum_{i\in I} a_i \vphi_i\, \right\|_{L_2}
\le C \ \mu \ (\log \mu)^{5/2}
\left\|\, \sum_{i\in I} a_i \vphi_i\, \right\|_{L_1}
\end{equation}
where  $\mu = L \, \sqrt {n/k} \,  \sqrt{\log k}$.

2) Moreover, if $n$ is an even integer, there exists a subset $I
\subset\{1, \ldots, n\}$ with $n/2 - c \sqrt n \le |I| \le  n/2 + c \sqrt n$
such that for every $a = (a_i) \in \C^n$,
$$
\left\|\, \sum_{i\in I} a_i \vphi_i\, \right\|_{L_2}
\le C \ L\,  \sqrt{\log n} (\log \log n)^{5/2}
\left\|\, \sum_{i\in I} a_i \vphi_i\, \right\|_{L_1}
$$
and
$$
\left\|\, \sum_{i\notin I} a_i \vphi_i\, \right\|_{L_2}
\le C \ L\,  \sqrt{\log n} (\log \log n)^{5/2}
\left\|\, \sum_{i\notin I} a_i \vphi_i\, \right\|_{L_1}.
$$
\end{theorem}

\noindent {\bf Remark.} Theorem \ref{Lambda_2} is satisfactory when
$k$ is an arbitrary proportion of $n$. It   strengthens   Theorem A
from \cite{GMPT} and extends the result of Talagrand
\cite{Tal_proport} which is applicable only when $|I|= n-k$ is a
sufficiently small proportion of $n$. Here we get  an arbitrary
proportion, with the distance $\sqrt {\log n} (\log \log n)^{5/2}$,
which is the right power of $\log n$ as we will see at the end of
the paper.

Define $\rho=\rho_{k,n}$, the restricted Kolmogorov $k$-width  of
the system   as the smallest number $\rho$  such that there exists a
subset $I \subset \{1, \dots, n\}$
   with $|I| \ge n - k $ satisfying
\begin{equation}
\left\|\, \sum_{i\in I} a_i \vphi_i\, \right\|_{L_2}
\le\rho
\left\|\, \sum_{i\in I} a_i \vphi_i\, \right\|_{L_1}
\end{equation}
for every $a = (a_i) \in \C^n$ (see \cite{GMPT} section 3). It was
proved in \cite {CT1} that $\rho=O(\sqrt{n/k}\log^3 n)$. This result
was improved  to $\rho=O(\sqrt{n/k}\sqrt{\log n}\,\log^{3/2} k)$ in
\cite{RudelsonVershynin}.

\bigskip

The proof of Theorem \ref{Lambda_2} is based on a random method
and  follows the argument given in \cite{GMPT} for proving Theorem
2.1. However,
instead of working with the space $L_1$, which is not uniformly
convex, we will approximate it by an $L_p$ space for $p$ "close" to
$1$ and use the full strength of the estimate given in Theorem
\ref{thrandom}.
\begin{proposition}
   \label{Lambda_2forp}
There exist two positive constants $c, C$ such that for
any orthonormal system $(\vphi_j)_{j=1}^n$ in $ L_2$ with
$\|\vphi_j\|_{L_\infty} \le L$ for $1 \le j \le n$, the following
holds.

1) For any $p \in (1,2)$ and any
   $1 < k < n$ there exists a subset $I \subset \{1, \dots, n\}$
   with $|I| \ge n - k $ such that for every $a = (a_i) \in \C^n$,
\begin{equation}
  \label{lambda_2_ineqforp}
\left\|\, \sum_{i\in I} a_i \vphi_i\, \right\|_{L_2} \le
\frac{C}{(p-1)^{5/2}} \, L \, \sqrt {n/k} \,  \sqrt{\log k} \left\|\,
\sum_{i\in I} a_i \vphi_i\, \right\|_{L_p}.
\end{equation}

2) Moreover, if $n$ is an even integer,
there exists a subset $I
\subset\{1, \ldots, n\}$ with $n/2 - c \sqrt n \le |I| \le  n/2 + c \sqrt n$
such that for every $a = (a_i) \in \C^n$,
$$
\left\|\, \sum_{i\in I} a_i \vphi_i\, \right\|_{L_2} \le
\frac{C}{(p-1)^{5/2}} \, L \, \sqrt {n/k} \,  \sqrt{\log k} \left\|\,
\sum_{i\in I} a_i \vphi_i\, \right\|_{L_p}
$$
and
$$
\left\|\, \sum_{i\notin I} a_i \vphi_i\, \right\|_{L_2} \le
\frac{C}{(p-1)^{5/2}} \, L \, \sqrt {n/k} \,  \sqrt{\log k} \left\|\,
\sum_{i\notin I} a_i \vphi_i\, \right\|_{L_p}.
$$
\end{proposition}

\medskip
\noindent{\bf Proof of Proposition~\ref{Lambda_2forp}.\ }
Let $X$ be the random vector taking the value $\vphi_i$ with
probability $1/n$ and denote by $E$ the complex vectorial space spanned
by $\{\vphi_1, \ldots, \vphi_n\}$.
Let $1 < k < n$ and let $X_1, \ldots, X_k$ be
independent copies of $X$. We define an operator
$\Gamma : E \to \ell_2^k$ by
$$
\forall y \in E, \quad
\Gamma y = \sum_{i=1}^k \lan X_i, y \ran e_i
$$
where $(e_1, \ldots, e_k)$ denotes the canonical basis of $\ell_2^k$.
Since $(\vphi_j)_{j=1}^n$ is an orthonormal system of $L_2$,
the basic properties of  $\Gamma$ are:

\smallskip
$\left\{
\begin{array}{ll}
    (i) & \dis \E \norm{\Gamma y}_{\ell_{2}^k}^2 = \frac{k}{n}  \,
    \|y\|_{L_2}^2,
    \\
    (ii) & \dis \ker \Gamma  = {\rm span} \{(\vphi_j)_{j=1}^n \setminus (X_i)_{i=1}^k \}.
\end{array}
\right.$

\medskip
We shall first prove that for any $\de \in (0,1)$, with probability
greater than $1- \de$, the set $(\vphi_{i})_{i \in I} =
\{(\vphi_j)_{j=1}^n \setminus (X_i)_{i=1}^k  \}$ satisfies
$$
\forall (a_i)_{i \in I} \in \C^{|I|},
\left\|\, \sum_{i\in I} a_i \vphi_i\, \right\|_{L_2} \le
\frac{C}{\de \, (p-1)^{5/2}} \, L \, \sqrt {n/k} \,  \sqrt{\log k} \left\|\,
\sum_{i\in I} a_i \vphi_i\, \right\|_{L_p}
$$
for a universal constant $C$.

Let $S_{L_2} = \{y : \|y\|_{L_2} = 1 \}$ be the unit sphere in $L_2$ and
observe that for any star-shaped subset $T \subset L_2$ the
following holds: if  $\rho > 0$ satisfies
\begin{equation}
  \label{rho}
  \sup_{y \in T\cap \rho S_{L_2}}
\left|\,   \sum_{j=1}^k \lan X_j, y\ran^2 \ -  \ \frac{k}{n} \, \rho^2 \, \right|
\le \frac{k \rho^2}{3 n},
\end{equation}
then
\begin{equation}
  \label{diameter}
\diam (\ker \Gamma \cap T) \le \rho.
\end{equation}
Indeed,
condition (\ref{rho}) implies that
for all $y \in T\cap \rho S_{L_2}$,
\begin{equation}
  \label{gamma}
\frac{2 k \rho^2}{3 n}
\le  \|\Gamma y \|_{\ell_2^k}^2
 \le \frac{4 k \rho^2}{3 n}.
\end{equation}

The homogeneity of (\ref{gamma}) and the fact that $T$ is star-shaped
imply  that if
the lower bound in
(\ref{gamma}) holds for all $y \in T\cap \rho S_{L_2}$,
then the same lower bound
 also holds for all $y \in T$ with $\|y\|_{L_2}\ge \rho$.
This in turn shows that if $y \in \ker \Gamma \cap T$ then
$\|y\|_{L_2} \le \rho$, as required in (\ref{diameter}).

Since $\ker \Gamma  = {\rm span} \{(\vphi_j)_{j=1}^n \setminus (X_i)_{i=1}^k \}$, we get
for
$T = B_{L_p} \cap E$
that whenever $\rho$ satisfies (\ref{rho}), then
$$
\left\|\, \sum_{i\in I} a_i \vphi_i\, \right\|_{L_2}
\le
\rho \left\|\, \sum_{i\in I} a_i \vphi_i\, \right\|_{L_p}
$$
for all scalars $a_i$, proving that \eqref{lambda_2_ineqforp} is
satisfied with the constant $\rho$.

\medskip

In order to find $\rho$ that satisfies \eqref{rho} with positive
probability we use Theorem \ref{thrandom}. Denote by $E_p$ the
complex Banach space spanned by $\vphi_1, \ldots, \vphi_n$ endowed
with the norm defined by
$$
\|y \| = \left( \frac{\|y\|_{L_p}^2 + \rho^{-2} \|y \|_{L_2}^2}{2} \right)^{1/2}.
$$
It is clear that $(B_{L_p} \cap E) \cap \rho B_{L_2} \subset
B_{E_p} \subset \sqrt
2 (B_{L_p} \cap E)  \cap \rho B_{L_2}$
and that the following properties are satisfied:

$\left\{
\begin{array}{l}
E_p \hbox{ is a Banach space with modulus of convexity of power
type 2}
\\
\hbox{with constant } \la^{-2} = p(p-1)/8,
\\
E_p^* \hbox{ is a Banach space of type 2 and } T_2(E_p^*) \le C \sqrt q =
C \sqrt{p/(p-1)}.
\end{array}
\right.
$

\noindent Indeed, the first property follows from Clarkson
inequality \cite{Clarkson}, that for any $f, g \in L_p$,
$$
\norm{ \frac{f+g}{2}}_{L_p}^2 + \frac{p(p-1)}{8} \norm{
\frac{f-g}{2}}_{L_p}^2 \le \frac{1}{2} (\|f\|_{L_p}^2 +
\|g\|_{L_p}^2).
$$
The second property is evident because for any $q \ge 2$, $L_q$ has
type 2 with constant $C \sqrt q$.
\\
By property $(i)$ of $\Gamma$, we get, taking $T = B_{L_p} \cap E$
and applying Theorem \ref{thrandom} to $E_p$,
\begin{eqnarray*}
\E  \sup_{y \in T\cap \rho S_{L_2}}
\left|\,   \sum_{j=1}^k  \lan X_j, y \ran^2  -  \frac{k}{n} \, \rho^2 \,
\right|
& = &
\E \sup_{y \in T\cap \rho S_{L_2}}
\left|\,  \sum_{j=1}^k  (\lan X_j, y \ran^2  -  \E \lan X_j, y \ran^2)  \,
\right|
\\
& \le &
 \E \sup_{y \in T\cap \rho B_{L_2}}
\left|\,  \sum_{j=1}^k  (\lan X_j, y \ran^2  -  \E \lan X_j, y \ran^2) \,
\right|
\\
& \le &
 \E  \sup_{y \in B_{E_p}}
\left|\,  \sum_{j=1}^k  (\lan X_j, y \ran^2  -  \E \lan X_j, y \ran^2) \,
\right|
\\
& \le & k (A^2 \, + \, A \, \sig)
\end{eqnarray*}
where $\sig  = \sup_{y \in B_{E_p}} \|y\|_{L_2} / \sqrt n \le \sqrt 2 \rho /
\sqrt n$ and
$$
\begin{array}{rcl}
A & = & \dis C \, \la^4  \, T_2(E_p^*) \sqrt{\frac{\log
k}{k}} \big(\E \max_{j \in J} \|\vphi_j\|_{E_p^*}^2\big)^{1/2}
\\
& \le & \dis C \, (p-1)^{-5/2} \ L \ \sqrt{\frac{\log k}{k}}
\end{array}
$$
since for every $1 \le j \le n$, $\|\vphi_j\|_{E_p^*} \le
\|\vphi_j\|_{L_{\infty}} \le L.$ By Chebychev inequality, for any $\de \in (0,1)$, we
conclude that
with probability greater than $1- \de$, there exists $X_1, \ldots, X_k$ such that
for any positive $\rho$,
$$
\sup_{y \in T\cap \rho S_{L_2}}
\left|\,   \|\Gamma y \|_{\ell_2^k}^2 \ - \ \frac{k}{n} \, \rho^2 \,
\right|
\le \frac{k}{\de}  (A^2 \, + \, A \, \sig).
$$
To conclude, we choose a constant $c$ large enough such that, for
$$
\rho = c  \sqrt{\frac{n}{k}} \ \sqrt{\log k} \ \frac{L}{\de \, (p-1)^{5/2}},
$$
the inequality $(\ref{rho})$ is satisfied with probability greater than $1-\de$.

The cardinality of the set $(\vphi_{i})_{i \in I} =
\{(\vphi_j)_{j=1}^n \setminus (X_i)_{i=1}^k  \}$ is greater than $n-k$ and
the first part of the proposition is proven choosing $\de = 1/2$.

\medskip
To prove the second part of the proposition, we make two more observations.
First, if  $\rho > 0$ satisfies $(\ref{rho})$ then we have proved that $(\ref{gamma})$
holds true. The upper estimate of this inequality implies that for all $y \in T \cap \rho S_{L_2}$,
$$
\sum_{i \in I} \lan \vphi_i, y \ran^2 \ge \sum_{i=1}^n \lan \vphi_i, y \ran^2 - \sum_{j=1}^k \lan X_j , y \ran^2
\ge \rho^2 \left(1 - \frac{4k}{3n} \right).
$$
Since $T$ is star shaped, we conclude as before that this inequality holds for all $y \in T$ for
which $\|y\|_{L_2} \ge \rho$. If $k < 3n/4$ then we have proved that if $y \in T$ and $\lan \vphi_i,y\ran = 0$
for all $i \in I$ (i.e. $y \in T \cap (\ker \Gamma)^{\perp}$), then
$\|y\|_{L_2} \le \rho$. But
$(\ker \Gamma)^{\perp} = {\rm span} \{(\vphi_i)_{i \notin I} \}$  and we
conclude that
if  $\rho > 0$ satisfies $(\ref{rho})$ with $k < 3n/4$ and $T = B_{L_p} \cap E$ then for any $(a_i)_{i=1}^n \in \C^n$,
$$
\left\|\, \sum_{i\in I} a_i \vphi_i\, \right\|_{L_2} \le
\rho
\left\|\,  \sum_{i\in I} a_i \vphi_i\, \right\|_{L_p}
$$
and
$$
\left\|\, \sum_{i\notin I} a_i \vphi_i\, \right\|_{L_2} \le
\rho
\left\|\,  \sum_{i\notin I} a_i \vphi_i\, \right\|_{L_p}.
$$
Secondly, it is not difficult to prove with a combinatorial argument (see Lemma 1.3 in \cite{GMPT})
that if $k = [\la n]$ with $\la = \log 2 < 3/4$ then with probability greater than $3/4$,
\begin{equation}
\label{cardinal}
n/2 - c \sqrt n \le |I| = n - |\{X_1, \ldots, X_k\}| \le n/2 + c \sqrt n,
\end{equation}
for some absolute constant $c >0$. Choosing $\de = 1/2$, we get that
with a positive probability,
both inequalities $(\ref{rho})$ and $(\ref{cardinal})$ are satisfied.
This concludes the
proof of the second point of the Proposition \ref{Lambda_2forp}.

\medskip

\noindent{\bf Proof of Theorem~\ref{Lambda_2}.\ } For any $p \in
(1,2)$, H\"older inequality states that for $\theta = (2-p)/p$,
$\|f\|_{L_p} \le \|f\|_{L_1}^{\theta} \|f\|_{L_2}^{1-\theta}$. Let
$\mu = L\, \sqrt {n/k} \,  \sqrt{\log k}$ and choose $p = 1 + 1 /
\log \mu$. Using Proposition \ref{Lambda_2forp}, there is a subset
$I$ of cardinality greater than $n-k$ for which
$$
\left\|\, \sum_{i\in I} a_i \vphi_i\, \right\|_{L_2} \le
C_p \, \mu \left\|\, \sum_{i\in I} a_i \vphi_i\,
\right\|_{L_p}
$$
where $C_p = C/(p-1)^{5/2}$. By the choice of $p$ and H\"older
inequality,
$$
\left\| \sum_{i\in I} a_i \vphi_i \right\|_{L_2}
\le
\mu \,
C_p^{p/(2-p)} \mu^{2(p-1)/(2-p)}
\left\|\, \sum_{i\in I} a_i \vphi_i \right\|_{L_1}
\le
C \, \mu \, h(\mu)
\left\| \sum_{i\in I} a_i \vphi_i  \right\|_{L_1}
$$
where $h(\mu) = (\log \mu)^{5/2}$ and $C$ is an absolute constant.

The same argument works for the second part of the Theorem~\ref{Lambda_2}.

\medskip
{\bf Appendix.}

As we mentioned in the Introduction, from a result of J. Bourgain
Theorem \ref{Lambda_2} is "almost" optimal.
We would like to thank
L. Rodriguez-Piazza for showing us the details of the proof of this optimality
in the case of the Walsh system.
We consider the Walsh system on $L_2[0,1]$ and take the $n$ first
functions
with $n = 2^N$.
In that case, it is well known that $(\vphi_1, \ldots, \vphi_n)$ form
a commutative multiplicative group that we denote by $(G, \cdot)$.
The main result is based on
constructing in any subset  of a group a translate of a subgroup of
big cardinality.

\begin{lemma}
    \label{subgroup}
Let $(G, \cdot)$ be the multiplicative group generated by the $2^N$
first Walsh functions.
For any $c \in (0,1)$ and for
any subset $\Lambda \subset G$ with $|\Lambda| \ge c \, 2^N$, we can
find  $b \in G$ and a subgroup $\Gamma$ of $G$ such that,
whenever $\log(1/c) \ge 1 /\,2^{N/2}$,

\smallskip
$\left\{
\begin{array}{ll}
    \medskip
   1. & \dis |\Gamma| = 2^p \ge N \log 2 \big/ \left( 3 \log(1/c)
   \right)
   \\
   2. & b \cdot \Gamma \subset \Lambda.
\end{array}
\right.
$
\end{lemma}

Assuming this result, we are able to prove the almost optimality of Theorem
\ref{Lambda_2} in the case of the Walsh system $(\vphi_1, \ldots,
\vphi_n)$ with $n = 2^N$.
Let $I$ be a subset of cardinality $n-k$ with $k \ge \sqrt n$. Taking $c = 1-
k/n$, we have  $\log(1/c) \ge k/n \ge 1/ \sqrt n$ and Lemma
\ref{subgroup} states that there exists $b \in (\phi_i)_{i =1}^n$ and a subgroup
$\Gamma$ of $G$ such that $b \cdot \Gamma \subset
(\phi_i)_{i \in I}$ and $|\Gamma| \ge (n \log n) / (20 k)$.
However, on any subgroup $\Gamma$ of  $G$, the $L_1$ norm and the
$L_2$ norm can not be compared with a better estimate than
$\sqrt{|\Gamma|}$ since
$$
\left\| \sum_{\gamma \in \Gamma}  \gamma \right \|_{L_2}
=
\sqrt{|\Gamma|}
\left\| \sum_{\gamma \in \Gamma}  \gamma \right \|_{L_1}.
$$
Therefore, inequality $(\ref{lambda_2_ineq})$ can not be satisfied
with a better constant than
$\sqrt{n \log n / 20 k}$. In particular when $k$ is proportional to
$n$,  the factor $\sqrt{\log n}$ is necessary in $(\ref{lambda_2_ineq})$.

\medskip
\noindent
{\bf Proof of Lemma \ref{subgroup}.}
We will prove the result in the case of the abelian additive group
$G=(\{0, 1\}^N, +)$. Note that this lands to a minor change of notation.
Let $n = 2^N$ be the cardinality
of this group.

Let $\gamma_1 \in G \setminus \{0\}$ be such that $|(\gamma_1 + \Lambda)
\cap \Lambda|$ is maximal and define $\Lambda_1 = (\gamma_1 + \Lambda)
\cap \Gamma$. Applying convolution, it is evident that
$$
\sum_{g \in G} |(g + \Lambda)  \cap \Lambda | = |\Lambda|^2
$$
and thus
$$
\begin{array}{rcl}
\dis |\Lambda_1| \ge
\frac{1}{n-1} \sum_{g \in G\setminus\{0\}}
|(g + \Lambda)  \cap \Lambda |
=
\frac{|\Lambda|^2 - |\Lambda|}{n-1}
& \ge &
\dis c^2 n \frac{1- 1/cn}{1-1/n}
\\
& \ge &
\dis c^2 n \left( 1 - \frac{1}{cn} \right).
\end{array}
$$
For notational convenience, we define $\gamma_0 = 0$ and $\Lambda_0 =
\Lambda$.
We iterate this procedure to construct a
family of points $\gamma_j$ and sets $\Lambda_j$ such that
$$
\left\{
\begin{array}{ll}
    (i) & \gamma_j \notin {\rm gr} \{\gamma_1, \ldots, \gamma_{j-1}\},
    \smallskip
    \\
    (ii) & \Lambda_j = (\gamma_j + \Lambda_{j-1})  \cap \Lambda_{j-1}
    \hbox{ has maximal cardinality,}
    \smallskip
    \\
    (iii) & \dis |\Lambda_j| \ge c^{2^j} n \left(1 - \frac{2^{j-1}}{n}
    \left(1/c + \ldots + 1/c^{2^{(j-1)}}\right) \right),
\end{array}
\right.
$$
where ${\rm gr} \{\gamma_1, \ldots, \gamma_{j-1}\}$ is the group
generated by $\{\gamma_1, \ldots, \gamma_{j-1}\}$. Once
$\gamma_{j-1}$ and $\Lambda_{j-1}$ are constructed, we again use the fact that
$$
\sum_{g \in G} |(g + \Lambda_{j-1})  \cap \Lambda_{j-1} | =
|\Lambda_{j-1}|^2
$$
and deduce that
$$
\sum_{g \in G\setminus {\rm gr}\{\gamma_1,
\ldots, \gamma_{j-1}\}}
|(g + \Lambda_{j-1})  \cap \Lambda_{j-1} | \ge |\Lambda_{j-1}|^2 -
2^{j-1} |\Lambda_{j-1}|.
$$
Therefore, there exists $\gamma_j \notin {\rm gr} \{\gamma_1,
\ldots, \gamma_{j-1}\}$ and $\Lambda_j = (\gamma_j + \Lambda_{j-1})  \cap \Lambda_{j-1}$
of maximal cardinality such that
$$
|\Lambda_j| \ge |\Lambda_{j-1}|
\frac{|\Lambda_{j-1}| - 2^{j-1}}{n-2^{j-1}}
$$
and $(iii)$ follows from a straightforward computation.

This construction can
be continued as long as
$$
1 - \frac{2^{p-1}}{n} \left( 1/c + \ldots + 1/c^{2^{(p-1)}} \right) >
0.
$$
This means that if $p$ is an integer such that
\begin{equation}
    \label{eq:group}
    n > p 2^{(p-1)} / c^{2^{(p-1)}}
\end{equation}
then $\Lambda_p \neq \emptyset$. Let $b \in \Lambda_p$, we shall prove
that
$$
b + {\rm gr}\{\gamma_1, \ldots, \gamma_p\} \subset \Lambda.
$$
Indeed let $x \in b + {\rm gr}\{\gamma_1, \ldots, \gamma_p\}$ and
define
$(\eps_1, \ldots, \eps_p) \in \{0,1\}^p$ such that
$x = b + \sum \eps_i \gamma_i$. Since $b \in \Lambda_p =
\Lambda_{p-1} \cap (\gamma_p + \Lambda_{p-1})$,  there exists
$\lambda_{p-1} \in \Lambda_{p-1}$ such that $b= \eps_p \gamma_p +
\lambda_{p-1}$. If one repeats the same argument for
$\lambda_{p-1}$ instead of $b$, then at the last step, we get $\lambda_1
\in \Lambda_1 = \Lambda \cap (\gamma_1 + \Lambda)$ and again, there
exists $\lambda \in \Lambda$ such that $\lambda_1 = \eps_1
\gamma_1 + \lambda$. Summarizing, we have found $\lambda \in \Lambda$
such that
$$
b = \lambda + \sum_{i=1}^p \eps_i \gamma_i.
$$
Therefore
$$
x = b + \sum_{i=1}^p \eps_i \gamma_i = \lambda \in \Lambda.
$$
To conclude, we  have to notice that if $p$ is chosen to be the
integer such that $2^p \ge \frac{N \log 2}{3 \log(1/c)} > 2^{p-1}$
then
$$
p 2^{(p-1)} \left(\frac{1}{c}\right)^{2^{(p-1)}}
< \frac{N^2 \log 2}{3 \log(1/c)}   \left(\frac{1}{c}\right)^{N \log 2
/ 3 \log(1/c)} \le N^2 2^{N/2 + N/3} \le 2^N
$$
using the fact that $\log(1/c) \ge 1 /\,2^{N/2}$.
The inequality $(\ref{eq:group})$ is
satisfied and  Lemma \ref{subgroup} holds true.

\address


\footnotesize
\begin{thebibliography}{Zz 99}
%
 \bibitem{BPST} J. Bourgain, A. Pajor, S.J. Szarek, N.
Tomczak-Jaegermann. {\it On the duality problem for entropy numbers of
operators.} Geometric Aspects of Functional Analysis (1987--88),
Springer Lecture Notes in Math. 1376 (1989), 50--63.
%
%
\bibitem{CT1} {E. Candes, T. Tao.} {\it Near optimal
    recovery from random projections: universal encoding strategies},
IEEE Inf. Theory 52 (2006), 5406-5425.

\bibitem{C} B. Carl. {\it Inequalities of Bernstein-Jackson-type and the degree of
 compactness of operators in Banach spaces.}
 Ann. Inst. Fourier (Grenoble)  {\bf 35}  (1985),  no. 3, 79--118.
%
\bibitem{Clarkson} J. A. Clarkson.
{\it Uniformly convex spaces,} Trans. Amer. Math. Soc.
{\bf 40} (1936), no.~3, 396--414
%
\bibitem{GineZinn}  E. Gin\'e,  J. Zinn. {\it Some limit
    theorems for empirical processes,} Ann. Probab.
  {\bf 12} (1984), no.~4, 929--998.
%
\bibitem{GMPT} O. Gu\'{e}don, S. Mendelson, A. Pajor, N.
Tomczak-Jaegermann. {\it Subspaces and orthogonal decompositions
  generated by bounded
  orthogonal systems.} Positivity, {\bf 11} (2007), no.~2, 269--283.

\bibitem{GueRud} O. Gu\'{e}don, M. Rudelson. {\it $L_p$
moments of random vectors via majorizing measures.} Adv. Math, {\bf 208}
(2007), no.~2, 798--823.
%
\bibitem{LinTza} J. Lindenstrauss,  L. Tzafriri. Classical Banach spaces. II.
Function spaces. Ergebnisse der Mathematik und ihrer Grenzgebiete [Results in
 Mathematics and Related Areas], 97. Springer-Verlag, Berlin-New
 York,  1979.
%
 \bibitem{Maurey} B. Maurey. {\it Type, cotype and $K$-convexity.}
 Handbook of the geometry of Banach spaces, Vol. 2,
 1299--1332, North-Holland, Amsterdam,  2003.
%
\bibitem{MaureyPisier} B. Maurey, G. Pisier.
{\it S\'eries de variables
al\'eatoires vectorielles ind\'ependantes et
propri\'et\'es g\'eom\'etriques des espaces de Banach.}
(French)  Studia Math.  58  (1976),  no. 1, 45--90.
%
 \bibitem{Pisier} G. Pisier.
    {\it Martingales with values in uniformly convex spaces,}
Israel J. Math. {\bf 20} (1975), no.~3-4, 326--350.
%
\bibitem{PisierBook} G. Pisier. The Volume of Convex Bodies and Banach
Space Geometry, Cambridge Tracts in Mathematics, 94, Cambridge
University Press, Cambridge, 1989.
%
\bibitem{Rod} L. Rodriguez-Piazza, {\it Personal communication}.
%
\bibitem{Rudelson} M. Rudelson.
{\it Random vectors in isotropic position}, MSRI preprint (1996).
%
\bibitem{Rudelson1} M. Rudelson.
{\it Random vectors in isotropic position,} J. Funct. Anal. {\bf
164} (1999), no.~1, 60--72.
%
\bibitem{RudelsonVershynin} M. Rudelson, R.Vershynin. {\it Sparse
reconstruction by convex relaxation: Fourier and Gaussian
measurements.} CISS 2006 (40th Annual Conference on Information
Sciences and Systems).
%
\bibitem{Sudakov} V.N. Sudakov, {\it Gaussian random processes and
measures of solide angles in Hilbert spaces}, Soviet. Math. Dokl.
{\bf 12} (1971), 412--415.
%
\bibitem{Talagrand} M. Talagrand. {\it
Sections of smooth convex bodies via majorizing
measures.}  Acta Math.  {\bf 175}  (1995),  no. 2, 273--300.
%
\bibitem{Tal:genericchaining} M. Talagrand.
{\it  Majorizing measures: the generic chaining.}
Ann. Probab. {\bf 24} (1996), no.~3, 1049--1103.
%
\bibitem{Tal_proport} M. Talagrand. {\it Selecting a proportion of
    characters}, Israel J. Math. {\bf 108} (1998), 173--191.
\end{thebibliography}
\end{document}